\documentclass{amsart}
\usepackage{geometry} 
\usepackage{amssymb}
\usepackage{cite}
\usepackage{graphicx}
\usepackage{url}
\usepackage{amsmath}
\usepackage{hyperref}
\usepackage{tikz}
\usetikzlibrary{matrix,arrows,positioning,automata,shapes}
\usepackage{verbatimbox}
\begin{document}

\newcommand{\Z}{\mathbb{Z}}
\newcommand{\mB}[1]{\mathbf{#1}}
\newcommand{\Bwr}{\boldsymbol{\wr}}

\newcommand{\SgpDec}{\textsc{SgpDec}}
\newcommand{\GAP}{\textsc{Gap}}
\newcommand{\Viz}{\textsc{Viz}}
\newcommand{\GraphViz}{\textsc{GraphViz}}

\pgfdeclarelayer{background layer}
\pgfsetlayers{background layer,main}

\newcommand{\lp}{\left(}
\newcommand{\rp}{\right)}

\title[SgpDec]{SgpDec: Cascade (De)Compositions of Finite Transformation Semigroups and Permutation Groups}
\author{Attila Egri-Nagy$^{1,3}$ \and James D. Mitchell$^{2}$ \and Chrystopher L. Nehaniv$^{3}$}
\address{
Centre for Research in Mathematics\\
School of Computing, Engineering and Mathematics\\
University of Western Sydney, Australia\\
\and
 School of Mathematics and Statistics\\
University of St Andrews, United Kingdom\\
\and
Centre for Computer Science \& Informatics Research\\
        University of Hertfordshire, United Kingdom}
\email{A.Egri-Nagy@uws.edu.au,jdm3@st-and.ac.uk,C.L.Nehaniv@herts.ac.uk}
\maketitle

\begin{abstract}
We describe how the \SgpDec~computer algebra package can be used for composing and decomposing permutation groups and transformation semigroups hierarchically by directly constructing substructures of wreath products, the so called cascade products.
\end{abstract}

\section{Introduction}

Wreath products are widely used theoretical constructions in group and semigroup theory whenever one needs to build a composite structure with hierarchical relations between the building blocks.
However, from a computational and engineering perspective they are less useful since wreath products are subject to combinatorial explosions and we are often interested only in substructures of them.
Cascade products precisely build these substructures by defining the hierarchical connections explicitly. 
As input, given a group or a semigroup  with unknown internal structure, the goal of cascade decomposition algorithms is to come up with a list of simpler building blocks and put them together in a cascade product, which realizes in some sense the original group or semigroup.
Roughly speaking, for permutation groups, cascade product decompositions can be interpreted as putting the inner workings of the  Schreier-Sims algorithm (generalized to any subgroup chain) into an external product form, therefore one can build cascade products isomorphic to the group being decomposed.
For semigroups, Krohn-Rhodes decompositions~\cite{primedecomp65} can be computationally represented by cascade products of transformation semigroups.   

In this paper we describe how the \GAP~\cite{GAP4} package \SgpDec~\cite{sgpdec} implements cascade products and decomposition algorithms and we also give a few simple example computations.
This description of the package only focuses on the core functionality of the package.

A \emph{transformation} is a function $f:X\rightarrow X$ from a set to itself,
and a \emph{transformation semigroup} $(X,S)$ of degree $n$ is a collection $S$ of transformations of $X$ closed under function composition, $|X|=n$. In case $S$ is a group of permutations of $X$, we call $(X,S)$ a \emph{permutation group}.
Using automata theory terminology sometimes we call $X$ the \emph{state set}, often represented as a set of integers $\mB{n}=\{0,\ldots,n-1\}$.
We write $x^s$ to denote the new state resulting from applying a transformation $s\in S$ to a state $x\in X$. 
\section{Cascade Product by a Motivating Example}
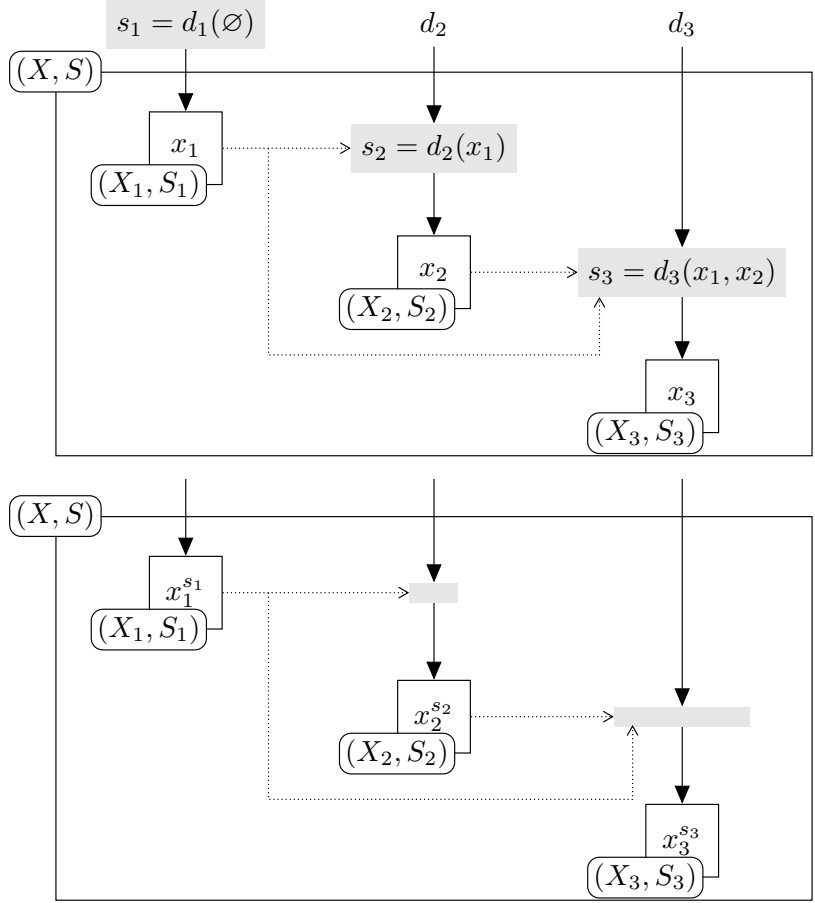
\begin{figure}
\begin{center}
\scalebox{1.1}{
\begin{tikzpicture}[node distance=40pt,->,>=triangle 45]
\tikzstyle{component} = [draw=black,rectangle,minimum size=25pt]
\tikzstyle{label} =[draw=black,rectangle, rounded corners,fill=white,inner sep=2pt, inner ysep=2pt]
\node[fill=black!10] (d1) at(0,4) {$s_1=d_1(\varnothing)$};
\node (d2) at (3,4) {$d_2$};
\node (d3) at (6,4) {$d_3$};

\node[fill=black!10] (g2) at (3,2.5) {$s_2=d_2(x_1)$};
\node[fill=black!10] (g3) at (6,1) {$s_3=d_3(x_1,x_2)$};

\node[component] (x1g1) at (0,2.5)  {$x_1$};
\node[label] at (x1g1.south west) {$(X_1,S_1)$};
\node[component] (x2g2) at (3,1) {$x_2$};
\node[label] at (x2g2.south west) {$(X_2,S_2)$};

\node[component] (x3g3) at (6,-0.5) {$x_3$};
\node[label] at (x3g3.south west) {$(X_3,S_3)$};
\node[component,inner sep=130pt,inner ysep=66pt] at (3,1.1) (xg) {};
\node[label] at (xg.north west) {$(X,S)$};
\tikzset{node distance=30pt};

\path[->]  (d1) edge (x1g1);
\path[>=angle 60,densely dotted]  (x1g1) edge (g2);
\path[->]  (d2) edge (g2);
\path[->]  (g2) edge (x2g2);
\path[->]  (g3) edge (x3g3);
\path[->]  (d3) edge (g3);
\path[>=angle 60,densely dotted]  (x2g2) edge (g3);
\draw[>=angle 60,densely dotted]  (1,2.5) -- (1,0) -- (5,0) -- (5,0.7);

\end{tikzpicture}}
\scalebox{1.1}{
\begin{tikzpicture}[node distance=40pt,->,>=triangle 45]
\tikzstyle{component} = [draw=black,rectangle,minimum size=25pt]
\tikzstyle{label} =[draw=black,rectangle, rounded corners,fill=white,inner sep=2pt, inner ysep=2pt]
\node (d1) at(0,4) {$\quad$};
\node (d2) at (3,4) {$$};
\node (d3) at (6,4) {$$};

\node[fill=black!10] (g2) at (3,2.5) {$\quad$};
\node[fill=black!10] (g3) at (6,1) {$\quad\quad\quad\quad$};

\node[component] (x1g1) at (0,2.5)  {$x_1^{s_1}$};
\node[label] at (x1g1.south west) {$(X_1,S_1)$};
\node[component] (x2g2) at (3,1) {$x_2^{s_2}$};
\node[label] at (x2g2.south west) {$(X_2,S_2)$};

\node[component] (x3g3) at (6,-0.5) {$x_3^{s_3}$};
\node[label] at (x3g3.south west) {$(X_3,S_3)$};
\node[component,inner sep=130pt,inner ysep=66pt] at (3,1.1) (xg) {};
\node[label] at (xg.north west) {$(X,S)$};
\tikzset{node distance=30pt};

\path[->]  (d1) edge (x1g1);
\path[>=angle 60,densely dotted]  (x1g1) edge (g2);
\path[->]  (d2) edge (g2);
\path[->]  (g2) edge (x2g2);
\path[->]  (g3) edge (x3g3);
\path[->]  (d3) edge (g3);
\path[>=angle 60,densely dotted]  (x2g2) edge (g3);
\draw[>=angle 60,densely dotted]  (1,2.5) -- (1,0) -- (5.4,0) -- (5.4,0.9);
\end{tikzpicture}}
\end{center}
\caption{Action in a cascade product of components $[(X_1,S_1)$, $(X_2,S_2)$, $(X_3,S_3)]$. The current state $(x_1,x_2,x_3)$ (top) is transformed to the new state $\lp x_1^{ s_1},x_2^{s_2},x_3^{s_3}\rp$ (bottom) by the transformation cascade $(d_1, d_2, d_3)$. The component actions $s_i$ are calculated by evaluating the dependency functions of $(d_1,d_2,d_3)$ on the states of the components above. The evaluations are highlighted and they happen at the same time. The dependencies, where the state information travels, are denoted by dotted lines.}
\label{fig:cascaction}
\end{figure}
To motivate the definition of the cascade product, we consider how the mod-4 counter, the cyclic permutation group $(\mB{4},\Z_4)$, can be constructed from two mod-2 counters.
The direct  product $\Z_2\times \Z_2$ contains no element of order 4. 
Since Aut$(\Z_2)$ is trivial there is only one semidirect product of $\Z_2$ and $\Z_2$, which equals their direct product. Their wreath product, $\Z_2\wr \Z_2 \cong D_4$,
the dihedral group of the square can be used to emulate a mod-4 counter, since $\Z_4\hookrightarrow D_4$.
But this construction is not efficient, beyond the required rotations the dihedral group has the flip-symmetry as well, doubling the size of the group. 
However, we would like to have a product construction that is isomorphic to  $(\mB{4},\Z_4)$.

This motivates the  definition of \emph{cascade products}: efficient constructions of substructures of wreath products, induced by explicit dependency functions~\cite{cascprod}. 
Essentially, cascade products are transformation semigroups glued together by functions in a hierarchical tree. 
More precisely, let $\big((X_{1},S_{1}),\ldots,(X_{n},S_{n})\big)$ be a fixed list of transformation semigroups, and dependency functions of the form
\[
d_i: X_1\times\ldots\times X_{i-1}\rightarrow S_i,\quad\text{for } i\in \{1,\ldots,n\}.
\]
A \emph{transformation cascade} is then defined to be an $n$-tuple of dependency functions $(d_1,\ldots,d_n)$, where $d_i$ is a dependency function of level $i$. If no confusion arises, on the top level we can simply write $d_1\in S_1$ instead of $d_1(\varnothing)\in S_1$. 
The cascade action is defined coordinatewise by  $x_i^{d_i(x_1,\ldots,x_{i-1})}$, applying the results of the evaluated dependency functions (see Fig.\ \ref{fig:cascaction}), so that the cascade product can be regarded as a special transformation representation on the set $X_1\times\ldots\times X_{n}$.
The hierarchical structure allows us to conveniently distribute computation among the components $(X_i, S_i)$, and perform abstractions and approximations of the system modelled as a cascade product.
Then if $W$ is a set of transformation cascades $(X_1,S_1)\wr_W\cdots \wr_W (X_n,S_n)$ denotes the transformation semigroup $(X_1 \times \cdots X_n, \langle W \rangle)$, where  $\langle W \rangle$ is the semigroup of transformation cascades generated by $W$.

\begin{figure}[t]
\begin{center}
\tikzset{->,>=triangle 45,auto,node distance=4cm}
\begin{tikzpicture}
\tikzstyle{every state}=[minimum size=3pt]
  \node[state]  (q_1)  {$0$};
  \node[fill=black!10]  (c0) [below left=0.7cm and 0.8cm of q_1] {$c(0)$};
  \node[state]  (q_2) [right of=q_1] {$1$};
  \node () [right=0.5cm of q_2] {top};
  \node[fill=black!10]  (c1) [below left=0.7cm and 0.5cm of q_2] {$c(1)$};
  \node[state]  (q_3) [below of=q_1] {$0$};
  \node[state]  (q_4) [right of=q_3] {$1$};
  \node () [right=0.5cm of q_4] {bottom};

  \path (q_1) edge  [bend left] node {$+1$} (q_2);
  \path (q_1) edge  [loop] node  [above] {$1_{\Z_2}$} (q_1);
  \path (q_2) edge  [bend left] node [above] {$+1$} (q_1);
  \path (q_2) edge  [loop] node  [above] {$1_{\Z_2}$} (q_2);
  \path (q_3) edge  [loop] node (id3) [above] {$1_{\Z_2}$} (q_3);
  \path (q_4) edge  [loop] node (id4) [above] {$1_{\Z_2}$} (q_4);
  \path (q_3) edge  [bend left] node (ct) {$+1$} (q_4);
  \path (q_4) edge  [bend left] node (cb) [above] {$+1$} (q_3);
  \path[>=angle 60,densely dotted]  (q_1) edge (c0);
  \path[>=angle 60,densely dotted]  (c0) edge (id3);
  \path[>=angle 60,densely dotted]  (c0) edge (id4);
  \path[>=angle 60,densely dotted]  (q_2) edge (c1);
  \path[>=angle 60,densely dotted]  (c1) edge (ct);
  \path[>=angle 60,densely dotted]  (c1) edge (cb);
\end{tikzpicture}
\caption{Two mod-2 counters cascaded together to build a mod-4 counter.}
\label{fig:mod4counter}
\end{center}
\end{figure}
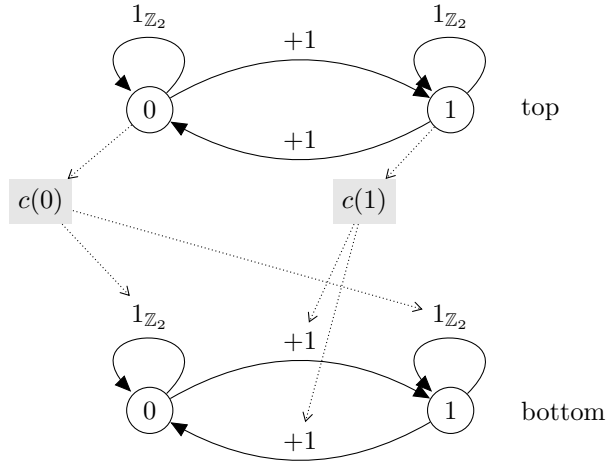

 We can construct $(\mB{4},\Z_4)$ exactly by using two copies of $(\mB{2},\Z_2)$. The generator set contains only one permutation cascade  $W=\{(+1,c)\}$, where $+1$ is the generator of $\Z_2$ and $c$ is a dependency function mapping $\mB{2}$ to $\Z_2$ with $c(0) =1_{\Z_2}$, and $c(1)=+1$.
The first dependency is a constant (increment modulo 2) while the second dependency implements the carry.
Therefore,  with fewer dependencies than required by the wreath product, the mod-4 counter can be realized by an isomorphic cascade product: $(\mB{2},\Z_2)\wr_W(\mB{2},\Z_2)\cong(\mB{4},\Z_4)$, see Fig.\ \ref{fig:mod4counter}.

An immediate consequence of the generality of the cascade product is that
several well-known constructions are special cases of the cascade product, and
as such they are easy to implement.
Direct products consist of all $d=(d_1,\ldots, d_n)$ with each $d_i$ constant. Wreath products consist of all possible dependency functions. 
Direct, cascade, and wreath products constructions for transformation semigroups are now available in \SgpDec, and iterated wreath products for permutation groups also became a bit more convenient to define.

\section{Functionality}

There are two different basic ways of using the \SgpDec~package. Depending on whether the starting point is a complex structure or a set of (simple) building blocks, we can do \emph{decomposition} or \emph{composition}.

\subsection{Composition and Construction}
\label{sect:composition}
The questions we aim to answer by constructing cascade products can be of the following types.

\begin{enumerate}
\item What is the (semi)group generated by a given set of transformation cascades? 
\item What can be built from a given set of (simple) components?
\end{enumerate}

\begin{figure}
\includegraphics[width=.45\textwidth]{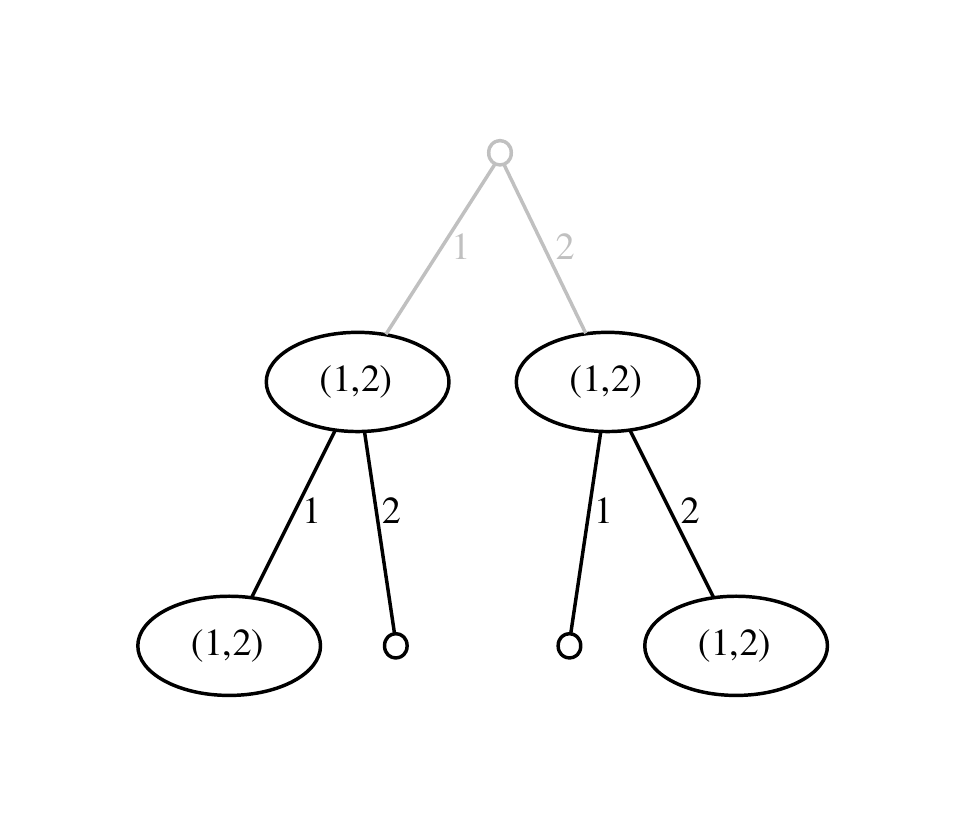}
\includegraphics[width=.45\textwidth]{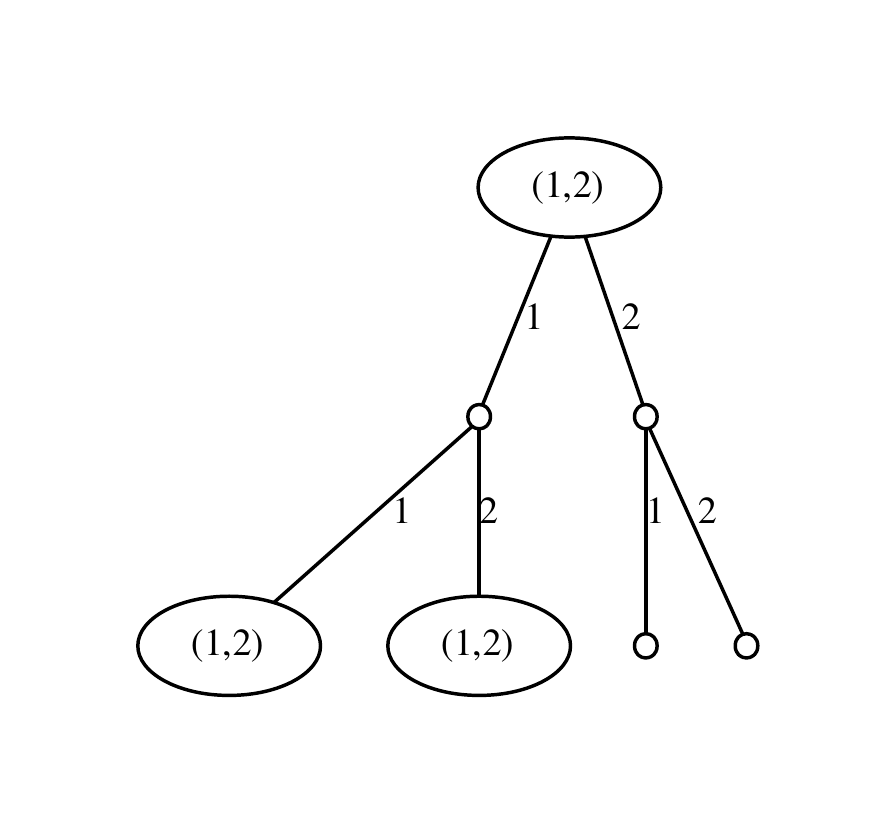}
\caption{Generators of a cascade representation of the quaternion group in a tree form. The edge labels are states, while the nodes contain the action. Empty node corresponds to the identity. The gray part of the tree is fixed.}
\label{fig:cascquaternion}
\end{figure}

The usual scenario is that for a list of components we give a set of cascades as a generating set.
For instance, the quaternion group $Q=\langle i , j \rangle$ is not a semidirect product, but it embeds into the full cascade product $(\mB{2},\Z_2)\Bwr(\mB{2},\Z_2)\Bwr(\mB{2},\Z_2)$, a group with 128 elements. Therefore, it can be built from copies of $\Z_2$.
The dependency functions can only have two values, thus to define cascade permutations it is enough to give only those arguments that give  $+1$ (the generator of $\Z_2$). A cascade permutation realizing $i$ is defined by the dependency functions $(d_1,d_2,d_3)$ where
$d_2(0)=d_2(1)=d_3(0,0)=d_3(1,1)=+1$
and all other arguments map to the identity. 
Similarly, a cascade realizing $j$ is defined by $(d'_1, d'_2,d'_3)$ where $d'_1(\varnothing)=d'_3(0,0)=d'_3(0,1)=+1$, (see Fig.\ \ref{fig:cascquaternion}, note that the state values are shifted by 1).
One can check that these two order 4 elements generate the 8-element quaternion group Q. Therefore by $W=\{(d_1,d_2,d_3),(d'_1,d'_2,d'_3)\}$ we have
$$(Q,Q)\cong (\mB{2},\Z_2)\wr_W(\mB{2},\Z_2)\wr_W(\mB{2},\Z_2).$$
\begin{verbbox}[\small]
gap> Z2:=CyclicGroup(IsPermGroup,2);
Group([ (1,2) ])
gap> d:=Cascade([Z2,Z2,Z2],[[[1],(1,2)],[[2],(1,2)],
                           [[1,1],(1,2)],[[2,2],(1,2)]]);
<perm cascade with 3 levels with (2, 2, 2) pts, 4 dependencies>
gap> dprime:=Cascade([Z2,Z2,Z2],[[[],(1,2)],[[1,1],(1,2)],[[1,2],(1,2)]]);
<perm cascade with 3 levels with (2, 2, 2) pts, 3 dependencies>
gap> StructureDescription(Group([d,dprime]));
"Q8"
\end{verbbox}
\theverbbox
\subsection{Decomposition}

\begin{enumerate}
\item What are the basic building blocks of a given (semi)group? 
\item How can we represent it as a cascade product?
\end{enumerate}

A typical scenario is that for a given composite semigroup or group we choose a decomposition algorithm which returns a cascade product. 
\subsubsection{Frobenius-Lagrange Decomposition.}
In the case of groups the decomposition uses the idea behind induction in representation theory (see e.g.\ \cite{AlperinBell}), so it traces back to Frobenius. Indeed, a special case of them comprises the well-known Krasner-Kaloujnine embeddings \cite{KrasnerKaloujnine}. All we need here is just standard group theory, namely the action on cosets, hence the name \emph{Frobenius-Lagrange Decomposition}.

How would someone come up with the generators cascades of the quaternion group in Section \ref{sect:composition}?
The easiest solution is to use this group decomposition.

\noindent\begin{verbbox}[\small]
gap> Q := QuaternionGroup(IsPermGroup,8);
Group([ (1,5,3,7)(2,8,4,6), (1,2,3,4)(5,6,7,8) ])
gap> CQ := FLCascadeGroup(Q);
<cascade group with 2 generators, 3 levels with (2, 2, 2) pts>
\end{verbbox}
\theverbbox

The actual implementation takes a subgroup chain as input (chief series by default) and form the components by examining the coset space actions derived from the chain.
Therefore, the decomposition method can be considered as a generalized Schreier-Sims algorithm~\cite{CGTHandbook}.

Coordinatewise calculation in a cascade product can also be thought of as a sequence of refining approximate solutions.
For instance, each completed step of an algorithm for solving the Rubik's Cube corresponds to calculating the desired value at a hierarchical level of some cascade product representation and it gives a configuration `closer' to the solved state.

\subsubsection{Holonomy Decomposition.}
\begin{figure}[th]
\begin{center}
\includegraphics[width=.75\textwidth]{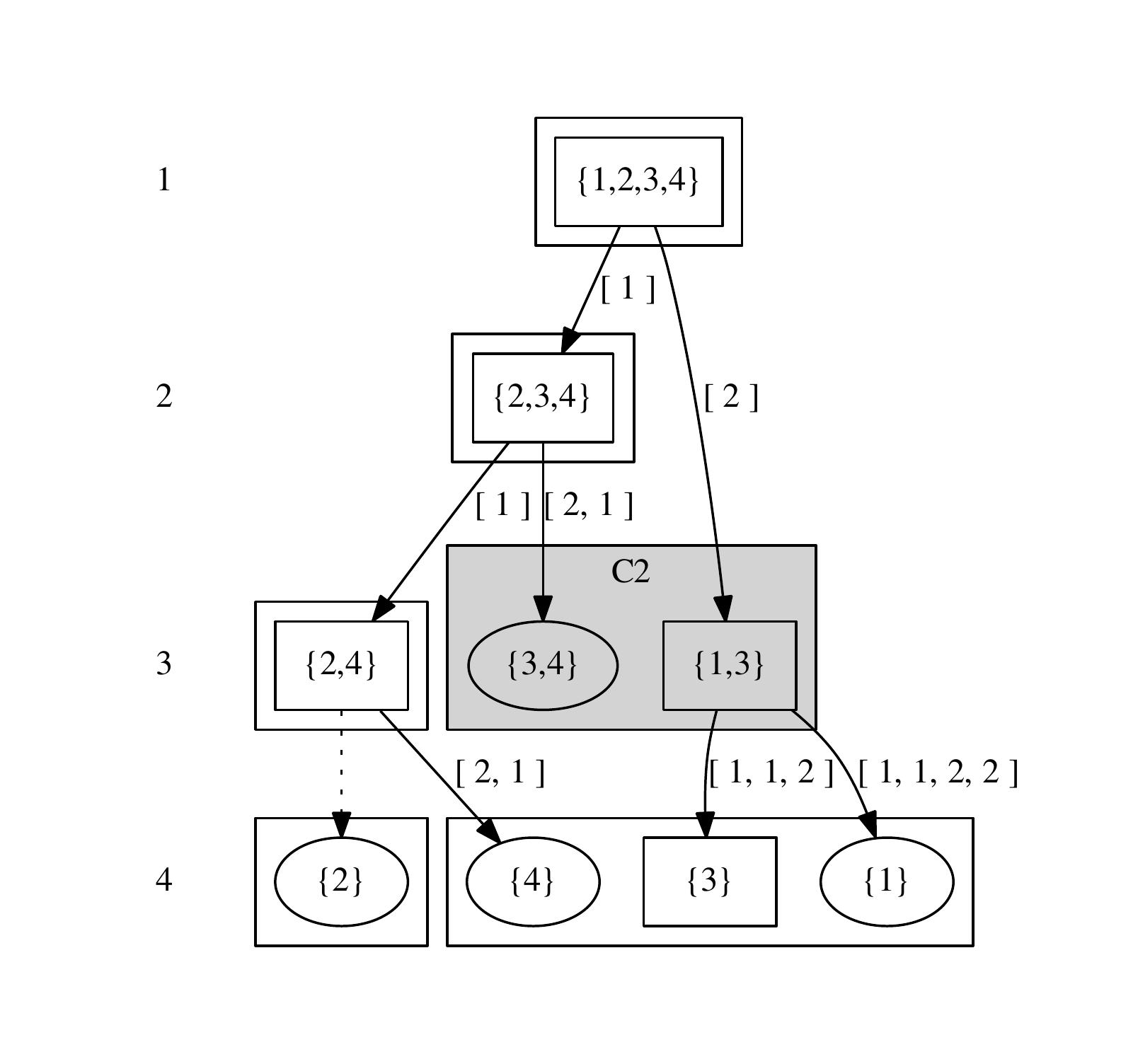}
\end{center}
\caption{Tiling picture -- the internal details of the holonomy decomposition of the transformation semigroup $T$ generated by $t_1$ and $t_2$. The numbers on the left denote the hierarchical levels (level 4 consists of singleton sets and it is needed by the holonomy algorithm but not a component of the cascade product). Outer boxes contain subsets that are mutually reachable from each other under the semigroup action. The arrows indicate how a subset is `tiled' by its subsets, the arrow labels contain words (sequences of generators) that take a subset to one of its tiles. Dotted arrow means the tile is not an image. Roughly, the holonomy algorithm finds the components by checking the action of the semigroup on a set of tiles.}
\label{fig:tiling}
\end{figure}
For transformation semigroups the holonomy method \cite{zeiger67a,zeiger68,ginzburg_book68,eilenberg,holcombe_textbook,KRTforCategories,automatanetworks2005} is  used.
The holonomy decomposition works by a close analysis of how the semigroup acts on those subsets of the state set which are images of the state set. 
As a small example let's define $T$ as the transformation semigroup generated by $t_1=\left(\begin{smallmatrix}1&2&3&4\\3& 2& 4&4\end{smallmatrix}\right)$ and $t_2=\left(\begin{smallmatrix}1&2&3&4\\3&3& 1& 3\end{smallmatrix}\right)$.
Calculating its holonomy decomposition and displaying some information can be done by the following commands:

\noindent\begin{verbbox}[\small]
gap> T:=Semigroup([Transformation([3,2,4,4]),Transformation([3,3,1,3])]);
<transformation semigroup on 4 pts with 2 generators>
gap> HT := HolonomyCascadeSemigroup(T);
<cascade semigroup with 2 generators, 3 levels with (2, 2, 4) pts>
gap> DisplayHolonomyComponents(SkeletonOf(HT));
1: 2 
2: 2 
3: (2,C2) 2 
\end{verbbox}
\theverbbox

\noindent The displayed information tells us that this 13-element semigroup can be realized as the cascade product of four copies of the transformation monoid of constant maps of two points and one instance of $\Z_2$. 
The components are put together in a 3-level cascade product.

Holonomy decompositions are useful whenever a finite state-transition model of some process needs to be analyzed (e.g.\ \cite{Dini2013}).

\subsection{Visualization}

\SgpDec~uses \GraphViz, a widely used graph drawing package \cite{Ellson03graphviz}, for visualisation purposes. 
The underlying idea is that a function generates  source code in the \texttt{dot} language for the given mathematical object.
Then the actual figure can be generated separately to be included in papers, or using the \Viz~package \cite{viz} immediately displayed on screen from the \GAP~command line.
Figure \ref{fig:cascquaternion} and \ref{fig:tiling} were both auto-generated using \GraphViz.
\section*{Acknowledgment}
The work reported in this article was funded in part by the EU project BIOMICS, contract number CNECT-318202. This support is gratefully acknowledged.

\bibliographystyle{plain}
\bibliography{coords}
\end{document}